\documentclass[12pt]{article}

\usepackage{mathrsfs,latexsym,amsmath}

\usepackage{pdfsync}

\usepackage{xypic}

\newcommand{\widebar}[1]{\overline #1}

\newcommand{\Pic}{\operatorname{Pic}}

\newcommand{\Q}{\mathbf Q}
\newcommand{\Z}{\mathbf Z}
\newcommand{\A}{\mathbf A}
\newcommand{\ind}{\operatorname{ind}}
\newcommand{\per}{\operatorname{per}}

\newcommand{\ms}[1]{\mathscr #1}

\newcommand{\m}{\bf \mu}

\newcommand{\Br}{\mathrm{Br}}

\newcommand{\Brp}{{}_p\mathrm{Br}}

\newcommand{\tensor}{\otimes}
\newcommand{\simto}{\stackrel{\sim}{\to}}

\title{Colliot-Th\'el\`ene's conjecture and finiteness of $u$-invariants}

\author{Max Lieblich \and R.\ Parimala \and  V. \ Suresh }
\date{}

\begin{document}
\maketitle

\section*{Introduction}

Let $X$ an excellent Noetherian regular scheme of dimension 2 which is
projective over an affine Dedekind scheme Spec$(A)$ and $K$ the
function field of $X$.  Let $l$ be a prime which is a unit on $X$
(i.e.  unequal to characteristic of the residue field of any point in
$X$). The work of Saltman on the division algebras over surfaces
([S1], [S3]) imply that given an element $\alpha$ in the $l$-torsion
of the Brauer group of $K$, there exist elements $f, g, \in K^*$ such
that $\alpha \otimes K(\sqrt[l]{f}, \sqrt[l]{g})$ is unramified on a
regular proper model of $K(\sqrt[l]{f}, \sqrt[l]{g})$ over $A$
(cf. [B]).  If $K$ is the function field of a curve over
a $p$-adic field and $l$ a prime not equal to $p$, the
unramified Brauer group on a regular proper model over $\mathbf{Z}_p$
is zero. Hence,  for $\alpha \in ~_lBr(K)$, index of $\alpha$ divides $l^2$.
 
In this paper, we  split the ramification of division algebras on
surfaces in a more general setting  without any assumption on the
characteristic  of the residue fields of points of the
scheme. More precisely, we prove the following (cf. Theorem 2.9): 

\paragraph*{Theorem.}Let $X$ be an excellent regular integral scheme
of dimension 2 and $K$ 
its function field. Suppose that char$(K) = 0$ and for every
codimension one point $x$ of $X$, if the characteristic of the residue
field $\kappa(x)$ at $x$ is $p$, then, $[\kappa(x) : \kappa(x)^p] = p$
(i.e. $p$-dimension of $\kappa(x)$ is 1).
If $\alpha$ is an element in the $p$-torsion of the Brauer group of
$K$, then there exists $f, g, h \in K^*$ such that $\alpha \otimes
K(\sqrt[p]{f}, \sqrt[p]{g}, \sqrt[p]{h})$ is unramified on a regular
model $Y$ of $K$ proper over $ X$. 

The above theorem leads to the following results on splitting ramification
of  exponent $p$ division algebras over function fields of curves over
$p$-adic fields and number fields (cf. Corollaries  2.10 and 2.11).

\paragraph*{Corollary.}If $k$ is a $p$-adic
field and $K$ a function field of a curve over $k$, then for every
element $\alpha \in ~_pBr(K)$, index$(\alpha)$ divides
period$(\alpha)^3$.

 \vskip 3mm

 \paragraph*{Corollary.}Let $k$ be a number field, $S$ the ring of
 integers in $k$, and $K$ the function field of a curve over $k$.  Let
 $\alpha \in Br(K)$ be a $p$-torsion element.There is an explicit
 degree $p^3$ extension $L$ of $K$ such that $\alpha \otimes L$ is the
 restriction of a class in the Brauer group of a regular proper model
 of $L$ over $S$.

\vskip 3mm

Thus the study of the period-index problem for such function fields is
reduced to a corresponding problem for unramified Brauer
classes. Recall the following conjecture of Colliot-Th\'el\`ene:

\vskip 3mm

\noindent 
{\bf CT-Conjecture}([CT1]): If $Y$ is a smooth projective
geometrically connected variety over a global field $k$ then the
Brauer-Manin obstruction to the existence of $0$-cycles of degree $1$
is the only obstruction.

 \vskip 2mm

 Using a moving lemma for $0$-cycles, we extend results of [L2] to
 show the following (cf. Theorem 3.3).

 \paragraph*{Theorem.} Let $k$ be a totally imaginary number field and $K$ the function
 field of a curve over $k$. Let $X$ be a regular proper model of $K$
 over the ring of integers in $k$.  If the CT-Conjecture holds then
 for any $\alpha \in {\rm Br}(X)$, the ind($\alpha$) divides
 per$(\alpha)^2$.

\vskip 2mm

Combining this with the results described above on splitting
ramification,  we deduce the following corollary (cf.  Theorem 4.1).

\paragraph*{Corollary.} Let $K$ be a function field in one variable
over a totally imaginary number field.  The CT-conjecture implies that for any $\alpha
\in Br(K)$, ind($\alpha$) divides period$(\alpha)^5$.  
\vskip 2mm

This together with a result on degree 3 Galois cohomology groups
([Su]) leads to the final result of the paper (cf. Theorem 4.5).

\paragraph*{Theorem.}If the CT-Conjecture holds then the function
field of a curve over a totally imaginary number field has finite
$u$-invariant.

\vskip 2mm

\noindent{\bf  Acknowledgements:}  We thank D. Saltman for generously sharing 
with us his notes on the Brauer group of complete discrete valuated fields in in equicharacteristic.
The first author is partially supported by the Sloan Foundation and National Science Foundation CAREER grant DMS-1056129.  
The second author is partially supported by National Science Foundation grant DMS-1001872. All authors are grateful to the American Institute of Mathematics, where this project was started as part of the workshop ``Deformation theory, patching, quadratic forms, and the Brauer group'' in January of 2012.

\section*{1. Brauer group of discretely valued fields}

Let $(K, \nu)$ be a complete discrete valued field with ring of
integers $R$, maximal ideal $m$ and residue field $\kappa$.  Suppose
that char$(K) = 0$, char$(\kappa) = p > 0$ and that $K$ contains a
primitive $p^{th}$ root of unity $\zeta$. Write
$N=\nu(p)p/(p-1)$. Since   
$$(p-1)\nu(\zeta-1)=\nu(p),$$   $N$ is a positive integer
divisible by $p$. 
Finally, for $i \geq 0$, let $$U_i = \{ u \in R^* \mid x \equiv 1 ~{\rm mod}~ m^i \}.$$
The following assumption about the residue field will play a key role throughout this paper.

\paragraph*{Assumption ($\ast$):} $[\kappa:\kappa^p] = p$. 

\paragraph{} Given $a, b \in K^*$, let $(a, b) \in \Brp(K)$ be the class of
the cyclic $K$-algebra defined by the relations
\begin{align*}
x^p &= a\\
y^p &= b\\
xy &= \zeta yx.
\end{align*}  
Let $br(K)_0 = \Brp(K)$  and for $i \geq
1$,  let $$br(K)_i\subset br(K)_0$$ be the subgroup generated by cyclic algebras
$(u, a)$ with $u \in U_i$ and $a \in K^*$.  Since $K$ is complete,
$br(K)_{n} = 0$ for $n > N$ ([CT2], 4.1.3  and [MS] ). In this section we
recall a few basic facts about $\Brp(K)$ and $br(K)_n$ due to Kato and
Saltman. 

\paragraph*{Lemma 1.1.}(Kato [K], Thm.\ 2) Let $$\alpha \in br(K)_m \setminus
br(K)_{m +1}$$ with $0 \leq m \leq N$.  Let $\pi \in R$ be a parameter
and $f \in R^*$ with $\overline{f} \not\in
\kappa^p$.
\begin{enumerate}
\item[(a)] If $m = 0$, then $$\alpha \equiv (u, \pi) \pmod{br(K)_1}$$ for some
$u \in R^*$ with  $\overline{u} \not\in \kappa^p$.
\item[(b)] If $m$ is coprime to $p$,  then there exists $x \in R^*$ such that  
$$\alpha \equiv (1 + \pi^m x, f) \pmod{br(K)_{m+1}}.$$
\item[(c)] If $0 < m < N$ and $m$ a multiple of $p$, then  $$\alpha \equiv (1 +
\pi^mx, \pi) \pmod{br(K)_{m+1}}$$ for some $x \in R^*$ with
$\overline{x} \not\in \kappa^p$.
\item[(d)] Let $b = \zeta-1$. 
If $m = N$ then $$\alpha = (1 + xb^p, f)  + ( 1 + x' b^p, \pi)$$ for some
$x,  x' \in R$.   Further $(1 + xb^p, f)$ is unramified at $\nu$.
\end{enumerate}

\paragraph*{Proof.}  Let $K_2(\kappa)$ be the Milnor $K$-group and
$$k_2(\kappa) = K_2(\kappa)/pK_2(\kappa).$$  The group $k_2(\kappa)$ is
isomorphic to a subgroup of $\Omega^2_\kappa$ (cf.\ [CT], 3.0). Since
(by Assumption ($\ast$)) we have $$[\kappa : \kappa^p] = p,$$ we know
that $$\Omega^2_\kappa = 0$$ and hence
$k_2(\kappa) = 0$.  Define a map
$$\kappa^*/\kappa^{*p} \to br(K)_0/br(K)_1$$ by $$\lambda \mapsto
(\tilde{\lambda}, \pi),$$ where for any $\lambda \in
\kappa$, $\tilde{\lambda} \in R$ is a lift. By ([K], Thm. 2, cf.\
[CT], 4.3.1), this map is an isomorphism. Hence every element in
$br(K)_0 \setminus br(K)_1$ is equivalent to $(u, \pi)$ modulo
$br(K)_1$ for some $u \in R^*$ with $\overline{u} \not\in \kappa^p$,
establishing part (a).

Let us prove part (b). Suppose $0 < m < N$ is coprime to $p$.
Define a map $$\Omega_{\kappa} \to br(K)_m/br(K)_{m+1}$$ by $$x\frac{dy}{y} \mapsto
(1 + \pi^m\tilde{ x}, \tilde{y}),$$
where $$x, y \in \kappa^*$$ are arbitrary elements with lifts 
$$\tilde{x}, \tilde{y} \in R^*.$$ By ([K], Thm. 2), this map is an isomorphism.
Since $\overline{f} \not\in \kappa^p$, we have $\kappa =
\kappa^p(\overline{f})$. In particular, since the (absolute)
differential vanishes on $\kappa^p$, every element in
$\Omega^1_\kappa$ has the form 
\begin{equation}\label{eq:1}
\omega = \sum_{i=1}^{p-1}\alpha_i\overline{f}^i\frac{d\overline{f}}{\overline{f}}
\end{equation}
for $\alpha_i\in\kappa^p$.
Thus every element of $br(K)_m/br(K)_{m+1}$ is the image of an element of the form $( 1 + x\pi^m,
f)$ for some $x \in R^*$.

Now we will prove part (c). Let $1 < m < N$ be divisible by $p$.
Define a map $$\kappa/\kappa^p \to
br(K)_m/br(K)_{m+1}$$ by $$\lambda + \kappa^p \to (1 + \pi^m
\tilde{\lambda}, \pi),$$
where for $\lambda \in \kappa$ we have chosen a lift $\tilde{\lambda}
\in R$.  Using \eqref{eq:1} above and the fact that $d$ is $\kappa^p$-linear, we see that every element of
$\Omega^1_\kappa$ is closed. Thus, the first summand in ([K], map
(iii) immediately preceding Thm.\ 2) is trivial. The map described
above is just the composition with the remaining summand of ([K], map
(iii) before Thm.\ 2) with the natural map $u_2^{(n)}/u_2^{(n+1)}\to
br_n(K)$ (in the notation of [K]).  Hence every non-zero
element in $br(K)_{m}/br(K)_{m+1}$ is the image of an element of the form $( 1 + \pi^mx, \pi)$
for some $x \in R^*$ with $\overline{x} \not\in \kappa^p$.

To prove part (d) (the case $n = N$),  
let $b=\zeta-1$. By ([K],
Thm.\ 2) or ([CT2], Thm.\ 4.3.1(d)) and using \eqref{eq:1} above, the element $\alpha$ can be
written in the desired form
$$\alpha = (1 + xb^p, f)  + ( 1 + x' b^p, \pi).$$ It remains to show
that the first summand is unramified (i.e., split by an unramified
extension of $K$). If $x\in\pi R$ then the first summand lies in
$br_{N+1}(K)=0$ and is unramified. 

Thus, we may assume that $x\in R^\ast$, so that $\nu(x)=0$.
We will show that in this case the field 
$$L=K(\sqrt[p]{1+xb^p})$$
is unramified over $K$; since this splits the algebra $(1+xb^p, f)$,
this will complete the proof.

Let $$\theta  = \sqrt[p]{1 + xb^p}$$ and $$d = \theta - 1.$$  
Let $\tilde{\nu}$ be the valuation on $L$ extending $\nu$ on
$K$. Since $1+xb^p$ is not a $p$th power, the minimal polynomial of
$\theta$ is $$g(z)=(z+1)^p-1-xb^p=0.$$
In particular, the norm of $\theta$ is $-xb^p$, so 
$$\tilde{\nu}(\theta)=\frac{1}{p}\nu(-xb^p)=\nu(b)$$
by our assumption that $x\in R^\ast$.
This means that $d = wb$ for some unit $w$ in the ring of integers in $L$.  It is easy to
see that $\overline{w}^p - \overline{w}  = \overline{x}$. Hence the $L/K$
is an unramified extension, as desired.  \hfill $\Box$

\paragraph*{Corollary 1.2.}  Every element $\alpha \in \Brp(K)$ is of
the form $(g, f) \cdot (h, \pi)$ for some $g, h \in R^*$.

\paragraph*{Proof.}  Let $0 \leq n \leq N$. Then, by (1.1), every
element in $br(K)_n$ is equivalent to $(g_n, f)(h_n, \pi)$ modulo
$br(K)_{n+1}$ for some $g_n, h_n \in R^*$. Since $br(K)_n = 0$ for $n
> N$, the corollary follows. \hfill $\Box$

\paragraph*{Lemma 1.3.} Let $n< N$ be coprime to $p$ and $c \in R^*$.
Then there exists $v \in R^*$ such that $(1 + \pi^nc, v\pi) = 1$.

\paragraph*{Proof.}  Since $n$ is coprime to $p$, there exists an
integer $m$ such that $nm \equiv 1 $ modulo $p$.  We have $$1 = (1 +
\pi^nc, -\pi^nc) ^m = (1 + \pi^nc, (-1)^m\pi^{nm}c^m) = (1 + \pi^nc,
v\pi)$$ with $v = (-1)^mc^m$. \hfill $\Box$

\paragraph*{Lemma 1.4.}  Let $n = mp<N$ and $i > 0$.  Let $u \in R^*$
be such that $\overline{u} \not\in \kappa^P$. Suppose $n+i $ is
coprime to $p$. Then every element in $br(K)_{n+i} \setminus
br(K)_{n+i+1}$ can be represented by $(1 + u\pi^n, u')$ for some $u'
\in R^*$.
 
\paragraph*{Proof.}  Fix $$\alpha\in br(K)_{n+1}\setminus
br(K)_{n+i+1}.$$
By (1.1) there is an $x\in R^\ast$ such that  $$\alpha \equiv (1 +
x\pi^{n+i}, u) \pmod{br(K)_{n+i+1}}.$$  Let $$ u' = (-1)^{p+1} (u - x\pi^i) ( 1 + x\pi^{n +
  i})^{-1}.$$ Since $i > 0$, $u \equiv u'$ modulo $\pi$ and hence $u'
= vu$ for some $v \in R^*$ with $v \equiv 1$ modulo $\pi$.  Thus, by ([CT],  4.1.1(b)), 
$$(1
+ x\pi^{n + i}, u) \equiv (1 + x\pi^{n + i}, u') \pmod{br(K)_{n+i
    +1}}.$$ 
Since $(1 - z^py, y) = 0$, we have $$(1 + x\pi^{n+i} , u' )
= ( (1 + x\pi^{n+i})(1- (-1)^p\pi^{mp}u'), u').$$ On the other hand,
we know $$(1 + x\pi^{n+i})(1 - (-1)^p\pi^{mp} u')= 1 + u\pi^{n}.$$
Thus $$\alpha \equiv (1 + u\pi^n, u') \pmod{br(K)_{n+i + 1}},$$ as
desired.  \hfill $\Box$

\paragraph*{Lemma 1.5.}  (Saltman [S4]) Let $\alpha \in br(K)_0
\setminus br(K)_1$.  Then there exist $u, v \in R^*$ such that
$\overline{u} \not\in \kappa^p$ and $$\alpha - (u, v\pi)$$ is
unramified.

\paragraph*{Proof.}  First we show by induction that for each $ 0 \leq
i < N$, there exist $u_i, v_i \in R^*$ such that $$\alpha - (u_i,
v_i\pi) \in br(K)_{i+1}$$ with $\overline{u}_i \not\in \kappa^p$.  By
(1.1), $$\alpha - (u_0, \pi) \in br(K)_1$$ for some $u_0 \in R^*$ with
$\overline{u}_0 \not\in \kappa^p$.  Suppose there exist $u_i, v_i \in
R^*$ such that $$\alpha - (u_i, v_i\pi) \in br(K)_{i+1}$$ and
$\overline{u}_i \not\in \kappa^p$.

Suppose that $i+1$ is coprime to $p$. Since 
$\overline{u}_i \not\in \kappa^p$, by (1.1), there exists
$c \in R^*$ such that $$\alpha - (u_i,  v_i\pi) - (u_i, 1 + \pi^{i+1}c) \in
br(K)_{i+2}.$$  Thus  $$\alpha -
(u_i, v_i\pi (1 + \pi^{i+1}c)) \in br(K)_{i+2}.$$ Suppose that $i+1$ is
divisible by $p$ and $i+1 < N$.  By (1.1),  there
exists $c \in R^*$ such that  $$\alpha - (u_i, v_i\pi) - ( 1 + \pi^{i+1}c,
v_i\pi) \in br(K)_{i+2}.$$   Then
$$\alpha - (u_i(1 + \pi^{i+1}c), v_i\pi) \in br(K)_{i+2},$$
and $\overline{u_i(1+\pi^{i+1}c)}=\overline{u}_i,$ so the induction
hypothesis is confirmed. 

In particular, $$\alpha - (u_{N-1}, v_{N-1}\pi) \in br(K)_{N}.$$ By
(1.1), $$\alpha - (u_{N-1}, v_{N-1}\pi) = \alpha' + (u, v_{N-1}\pi)$$
with $\alpha'$ unramified and $u \in 1+\pi R$.  Thus $$\alpha - (u_{N-1}u,
v_{N-1}\pi) = \alpha'$$ is unramified and $\overline{u_{N-1}u}\not\in\kappa^p$, as desired. \hfill $\Box$
  
\paragraph*{Lemma 1.6.}  (Saltman [S4]) Let $1 \leq n < N$ be coprime
to $p$ and $\alpha \in br(K)_n \setminus br(K)_{n+1}$.  Then there
exist $u, v \in R^*$ such that $$\alpha - (1 + v\pi^n, u\pi)$$ is
unramified.

\paragraph*{Proof.}  First we show by induction that for each $ 0 \leq
i < N- n$, there exist $u_i, v_i \in R^*$ such that $$\alpha - (1 +
v_i\pi^n, u_i ) \in br(K)_{n+i+1}$$ with $\overline{u}_i \not\in
\kappa^p$.  Since $n$ is coprime to $p$, by (1.1), $$\alpha - (1 +
v_0\pi^n, u_0) \in br(K)_{n+1}$$ for some $u_0, v_0 \in R^*$ with
$\overline{u}_0 \not\in \kappa^p$.  Suppose there exist $u_i, v_i \in
R^*$ such that $$\alpha - (1 + v_i\pi^n, u_i) \in br(K)_{n+ i+1}
\setminus br(K)_{n+i+2}$$ and $\overline{u}_i \not\in \kappa^p$.  We
will break the proof into two cases, depending upon the divisibility
properties of $n+i+1$.

\textbf{Case 1:} $n+i+1$ is coprime to $p$. Since
$\overline{u}_i \not\in \kappa^p$, by (1.1), there exists $c \in R^*$
such that $$\alpha - (1 + v_i\pi^n, u_i) - (1 + \pi^{n+i+1}c, u_i) \in
br(K)_{n+i+2}.$$ Let $$v_{i+1} = v_i + \pi^{i+1}c + \pi^{n+i+1}
v_ic.$$ Then $$\alpha - (1 + v_{i+1}\pi^n, u_i ) \in br(K)_{n+ i+2},$$
verifying the inductive hypothesis.

\textbf{Case 2:} $n+i+1$ is divisible by $p$ and $n+i+1 < N$. Since $n$ is
coprime to $p$, by (1.3) there is a $v\in R^\ast$ such that 
$$(1 + v_i\pi^n, v\pi)=1.$$  In particular, $$(1 + v_i \pi^n, u_i) = (1 + v_i \pi^n,
u_iv\pi).$$  By (1.1), there exists $c \in R^*$ such that $$\alpha - (1
+ v_i\pi^n, u_i) - ( 1 + \pi^{n+i+1}c, u_iv\pi) \in br(K)_{n+i+2}.$$
Thus 
\begin{align*}\alpha - (1 + v_i\pi^n, u_i) - ( 1 + \pi^{n+i+1}c, u_iv\pi) &=
\alpha - (1 + v_i\pi^n, u_iv\pi) - ( 1 + \pi^{n+i+1}c, u_iv\pi) \\
&=
\alpha - ( (1 + v_i\pi^n)(1 + \pi^{n+i+1}c), u_iv\pi)\\
 &\in br(K)_{n+
  i+2}.
\end{align*}

Let $v_{i+1} = v_i + \pi^{i+1}c + v_ic\pi^{n+i+1}$. Then
$$\alpha - (1+ v_{i+1}\pi^n, u_iv\pi) \in br(K)_{n+i+2}.$$ Once again,
by (1.3), $$(1 + v_{i+1}\pi^n, v'\pi) = 1$$ for some $v' \in R^*$.
Hence $$(1 + v_{i+1}\pi^n, v'^{-1}\pi^{-1}) = 1$$ and 
\begin{align*}
(1 +
v_{i+1}\pi^n, u_iv\pi) &= (1 + v_{i+1}\pi^n, u_iv\pi) (1 +
v_{i+1}\pi^n, v'^{-1}\pi^{-1}) \\
&= (1 + v_{i+1}\pi^n, u_ivv'^{-1}).
\end{align*}
Let $u_{i+1} = u_ivv'^{-1}$. Then $$\alpha - (1 + v_{i+1}\pi^n,
u_{i+1}) \in br(K)_{n+i+2}.$$ 
Suppose that $\overline{u}_{i+1} \in \kappa^p$.  
Replacing $u_{i+1}$ by $a^pu_{i+1}$, we assume that $\overline{u}_{i+1} = 1$.
Then, by ([CT], 4.1.1(b)), $(1 + v_{i+1}\pi^n, u_{i+1}) \in br(K)_{n+1}$, 
contradicting the fact that  $\alpha \not\in br(K)_{n+1}$. Thus
 $\overline{u}_{i+1} \not\in \kappa^p$,  thereby  verifying the
inductive hypothesis.

In particular, for $i=N-n$, we have $$\alpha - (1 + v_{N-1}\pi^n,
u_{N-1}) \in br(K)_N.$$ Since $n$ is coprime to $p$, by (1.3), we have
$$(1 + v_{N-1}\pi^n, u_{N-1}) = (1 + v_{N-1}\pi^n, u_{N-1}v\pi)$$ for
some $v \in R^*$. By (1.1(d)), 
$$\alpha - ( 1+v_{N-1}\pi^n, u_{N-1}v\pi)
= \alpha'  + (1 + x'b^p, u_{N-1}v\pi)$$ for some $\alpha'$ unramified
and $x' \in R$. Hence $$\alpha - ( (1+v_{N-1}\pi^n)(1 + x'b^p),
u_{N-1}v\pi) = \alpha'$$ is unramified, as desired.  \hfill $\Box$
  
\paragraph*{Lemma 1.7.} (Saltman [S4]) Let $1 < n < N$ be divisible by
$p$ and $$\alpha\in br(K)_n \setminus br(K)_{n+1}.$$  Then
there exist $u, v \in R^*$ such that $$\alpha - (1 + u\pi^n, v\pi)$$
is unramified.

\paragraph*{Proof.}  First we show by induction that for each $ 0 \leq
i < N- n$, there exist $u_i, v_i \in R^*$ such that $$\alpha - (1 +
u_i\pi^n, v_i\pi ) \in br(K)_{n+ i+1}$$ and $$\overline{u}_i \not\in
\kappa^p.$$ Since $n$ is divisible by $p$, by (1.1) there exists $ u_0 \in R^*$ with
$\overline{u}_0 \not\in \kappa^p$ such that $$\alpha - (1 +
u_0\pi^n, \pi) \in br(K)_{n+1}.$$  Suppose there exist $u_i, v_i \in
R^*$ such that $$\alpha - (1 + u_i\pi^n, v_i\pi) \in br(K)_{n+ i+1}
\setminus br(K)_{n+i+2}$$ and $\overline{u}_i \not\in \kappa^p$.  We
again break into two cases.

\textbf{Case 1:} $n+i+1$ is divisible by $p$. By (1.1), there exists
$c \in R^*$ such that $$\alpha - (1 + u_i\pi^n, v_i\pi) - (1 +
\pi^{n+i+1}c, v_i\pi) \in br(K)_{n+i+2}.$$ Thus
\begin{align*}
  \alpha - ((1 + u_i\pi^n)(1 + \pi^{n+i+1}c), v_i\pi ) &= \alpha - (1
  + (u_i +
  \pi^{i+1}c + u_i\pi^{n+i+1}c)\pi^n, v_i\pi)\\
  &\in br(K)_{n+ i+2}.
\end{align*} 
Let
$$u_{i+1} = u_i + \pi^{i+1}c + u
_i\pi^{n+i+1}c.$$ Then $\overline{u}_{i+1} = \overline{u}_i \not\in
\kappa^p$ and $$\alpha - (1 + u_{i+1}\pi^n, v_i\pi) \in
br(K)_{n+i+2},$$ establishing the inductive hypothesis.

\textbf{Case 2:} $n+i+1$ is coprime to $p$.  By (1.1), there are $c,f
\in R^*$ with $\overline{f} \not\in \kappa^p$ such that $$\alpha - (1 +
u_i\pi^n, v_i\pi) - ( 1 + c\pi^{n+i+1}, f) \in br(K)_{n+i+2}.$$ By
(1.4), there is an $f'\in R^\ast$ such that $$(1 +
c\pi^{n+i+1}, f) = ( 1 + v_i\pi^n, f').$$  In
particular 
\begin{align*}
\alpha - (1 + u_i\pi^n, v_if'\pi) &=
\alpha -(1 +u_i\pi^n, v_i\pi) - (1 + u_i\pi^n, f') \\
&= \alpha - (1+ u_i\pi^n, v_i\pi) - (1 +c\pi^{n+i+1}, f) \\
 &\in br(K)_{n+i+2}.
\end{align*}
 
Proceding as in the proof of (1.5), we conclude that there are
$u,v\in R^\ast$ such that $$\alpha -
(1 + u\pi^n, v\pi)$$ is unramified. \hfill $\Box$

\paragraph*{Proposition 1.8.} (Saltman [S4]) Let $\alpha \in
\Brp(K)$. Then there exists a parameter $\pi'$ such that $\alpha \otimes
K(\sqrt[p]{\pi'})$ is unramified.

\paragraph*{Proof.}  If $\alpha\in br(K)_n$ for $0\leq n < N$ $\alpha
\in br(K)_N$, then applying (1.5, 1.6, 1.7) we see that there is a
uniformizer $\pi$ and an element $a\in K$ such that 
$$\alpha-(a,\pi)$$
is unramified. If $\alpha\in br(K)_N$ then by (1.1) we have (in fact) that there
is an element $x\in R$ such that 
$$\alpha - (1+x\pi^N,\pi)$$
is unramified for any choice of uniformizer $\pi$.  In any case, this
implies that for the appropriate choice of $\pi$, the base change 
$$\alpha \otimes K(\sqrt[p]{\pi})$$ is unramified, as desired.
\hfill $\Box$.

\section*{2. The Brauer group of the  function field of a surface}  

Let $R$ be an integral domain with field of fractions $K$. 
Let $C$ be a central simple algebra over $K$. 
\paragraph*{Definitions.} 
\begin{itemize}
\item The algebra 
$C$ is {\it unramified} on $R$ if there is an Azumaya algebra ${\cal
  C}$ over $R$ such that ${\cal C} \otimes_R K$ is Brauer equivalent
to $C$.  
\item An element in $\Br(K)$ is {\it unramified} on $R$
if it is represented by central simple algebra over $K$ which is
unramified on $R$.  
\item If an element in $\Br(K)$ is not unramified on
$R$, then we say that it is {\it ramified} on $R$.  
\item Let $P$ be a prime
ideal of $R$. If an element in $\Br(K)$ is unramified (resp.\ ramified)
on $R_P$, then we say that $\alpha$ is unramified (resp.\ ramified) at
$P$.
\end{itemize}

\paragraph{}
Let $X$ be a regular integral scheme and $K$ its function field.
Let $X^i$ be the set of points of $X$ of codimension $i$. For $x \in
X^i$, let ${\cal O}_{X, x}$ be the local ring at $x$, $m_x$ be the
maximal ideal at $x$ and $\kappa(x)$ be the residue field at $x$.  For
$x \in X^1$, if a central simple algebra $C$ over $K$ is unramified
(resp.\ ramified) at ${\cal O}_{X, x}$, then we say that $C$ is
unramified (resp.\ ramified) at $x$.  Let $K_x$ denote the field of
fractions of the completion of ${\cal O}_{X, x}$ at $m_x$.  If $m_x$
is generated by $\pi$, we also denote $K_x$ by $K_{\pi}$.  For $x \in
X^1$, let $\nu_x$ be the valuation on $K$ given by $x$.  Suppose that
char$(K) = 0$ and char$(\kappa(x))= p > 0$.  Let $N_x =
\nu_x(p)p/(p-1)$.  If $m_x$ is generated by $\pi$, we also denote
$N_x$ by $N_{\pi}$.

\vskip 5mm 
We will make (implicit) essential use of the following results
throughout the rest of this paper. The following results are true in much more generality 
for  torsors (cf. [CTS], 6.13, [APS], 4.3).

\paragraph*{Lemma 2.1}
Suppose $X$ is a Noetherian two dimensional regular integral scheme
with function field $K$. Given $\alpha\in\Br(K)$, the following are
equivalent.
\begin{enumerate}
\item For every discrete valuation of $K$ with center on $X$ and
  valuation ring $R$, $\alpha$
  is unramified on $R$.
\item For every discrete valuation of $K$ with center of codimension
  $1$ on $X$ and valuation ring $R$, $\alpha$ is unramified on $R$.
\item For any central simple algebra $C$ over $K$ representing
  $\alpha$, there is a sheaf of Azumaya algebras $A$ on $X$ and an
  isomorphism $A\tensor_{\ms O_X}K\simto C$ of $K$-algebras.
\end{enumerate}
\paragraph{Proof.}
Fix a central simple algebra $C$ over $K$ representing $\alpha$ and
let $A$ be a maximal order in $C$. By [Prop.\ 1.2, AG] any localization of $A$ at a point of $X$ is a maximal
order in $C$ over $\ms O_{X,x}$. 
Since $X$ is regular and 2-dimensional and $A$ is reflexive 
as an $R$-module [Thm.\ 1.5, AG], $A$ is
a finite locally free $R$-algebra. Consider the map
$$\mu:A\otimes A^\circ \to\ms End_{\ms O_X}(A)$$
of locally free sheaves of $\ms O_X$-modules of rank $n^2$ associated
to the map of presheaves that sends an elementary
tensor $a\otimes b$ in $A(U)\tensor_{\ms O(U)}A(U)$ to the $\ms 
O_X(U)$-linear map sending $x$ to $axb$. The maximal order $A$ is Azumaya if
and only if $\mu$ is an isomorphism.
We therefore assume that $X = {\rm Spec}(B)$ with $B$ a regular local ring of dimension 2. 
In this case $A \tensor_BA^\circ $ and $End_B(A)$ are free $B$-modules of same rank, since
$\mu \tensor_B K$ is an isomorphism. For a choice of a basis of these modules over $B$, $\mu$ 
is an isomorphism if and only if the determinant of $\mu$ is a unit in $B$. 
By Krull's theorem, the determinant is a unit  in $B$ if and only if each
it is a  unit at every codimension one point. This shows that the second two
conditions are equivalent.

On the other hand, the first condition certainly implies the second,
and hence the third. It remains to show that the third condition
implies the first. Suppose $A$ is an Azumaya algebra on $X$
restricting to $C$. If $R$ is a valuation ring with center on $X$ then
there is a commutative diagram
$$\xymatrix{& \operatorname{Spec} K\ar[dr]\ar[dl] & \\
\operatorname{Spec} R\ar[rr] & & X 
}$$ 
in which the two diagonal arrows are the canonical ones. Pulling $A$
back yields an Azumaya algebra on $R$ restricting to $C$ on $K$,
showing that $C$ is unramified on $R$, as desired.
\hfill $\Box$
\vskip 5mm

\paragraph*{Corollary 2.2}
Suppose that $R$ is a two dimensional regular Noetherian integral
domain with field of fractions $K$.    Then a central simple algebra 
$C$ over $K$ is unramified on $R$ if and only if it is unramified on 
$R_P$ for every height one prime ideal $P$ of $R$.

\paragraph*{Corollary 2.3}
Suppose that $X$ is a Noetherian two dimensional regular scheme with
fraction field $K$. Then an element $\alpha\in\Br(K)$ is unramified if
and only if it is in the image of the injective restriction map
$$\Br(X)\to\Br(K).$$
\vskip 5mm

\paragraph*{Lemma 2.4} Let $R$ be a discrete valuation ring with field
of fractions $K$. Let $\hat{R}$ be the completion of $R$ at the
discrete valuation and $\hat{K}$ the field of fractions of $\hat{R}$.
Then a central simple algebra $C$ over $K$ is unramified at $R$ if and
only if $C \otimes_K \hat{K}$ is unramified at $\hat{R}$.

\paragraph*{Proof.}  Let ${\cal A}$ be an Azumaya algebra over
$\hat{R}$ such that there is an isomorphism $\phi: {\cal A}
\otimes_{\hat{R}} \hat{K} \to C \otimes_K \hat{K}$.  Then ${\cal C} =
{\cal A} \times _\phi C = \{ (a, b) \in {\cal A} \times C \mid \phi (a
\times 1) = (b, 1) \}$ is an Azumaya $R$-algebra with ${\cal C}
\otimes_RK \simeq C$.  \hfill $\Box$

\paragraph*{Lemma 2.5} Let $R$ be a discrete valuation ring with field
of fractions $K$ and residue field $F$. Suppose that char$(F) \neq p$
and $\alpha \in \Brp(K)$. If $\pi$ is a parameter of $R$, $\alpha
\otimes_KK(\sqrt[p]{\pi})$ is unramified on $R[\sqrt[p]{\pi}]$.

\paragraph*{Proof.}    See ([S2], 0.4). \hfill $\Box$

\paragraph*{Proposition 2.6.}  Let $A$ be a regular local ring of
dimension two with maximal ideal $m = (\pi, \delta)$. Let $K$ be the
field of fractions of $A$ and $k$ the residue field of $A$.  Suppose
that char$(K) = 0$, char$(k) = p>0 $ and $k = k^p$. 
Let $\alpha \in \Brp(K)$. If $\alpha$ is
ramified on $A$ at most at $(\pi)$ and $(\delta)$, then $\alpha
\otimes_K K(\sqrt[p]{\pi}, \sqrt[p]{\delta})$ 
is unramified at any discrete valuation  of $K(\sqrt[p]{\pi},
\sqrt[p]{\delta})$ dominating $A$. 

\paragraph{Proof.}  Let $\nu$ be a discrete valuation of $L =
K(\sqrt[p]{\pi}, \sqrt[p]{\delta})$ dominating $A$. Let $R$ be the
valuation ring at $\nu$.  The completion of $A$ at $m$ is contained in
the completion of $R$ at its maximal ideal. By (2.4), to show that
$\alpha$ is unramified at $R$, we replace $A$ by its completion at $m$
and assume that $A$ is complete.

Let $B = A[\sqrt[p]{\pi}, \sqrt[p]{\delta}]$. Then $B$ is a complete
regular local ring with field of fractions $L$ and $B$ is the integral
closure of $A$ in $L$. Hence $B \subset R$.  Thus, to show that
$\alpha$ is unramified at $\nu$, it is enough to show that $\alpha$ is
unramified on $B$.  Since $B$ is a regular ring of dimension 2, by
(2.2), it is enough to show that
$\alpha$ is unramified at every height one prime ideal of $B$.

Let $Q$ be a height one prime ideal of $B$.  Since $B$ is integral
over $A$, $Q \cap A = P$ is a height one prime ideal of $A$. Suppose
$P \neq (\pi)$ and $P \neq (\delta)$.  Then $\alpha$ is unramified at
$P$ and hence $\alpha$ is unramified at $Q$.

Suppose $P = (\pi)$.  Then $L_Q = K_P(\sqrt[p]{\pi},
\sqrt[p]{\delta})$.  Suppose that char$(\kappa(P)) \neq p$.  Since
$\pi$ is a parameter at $P$, by (2.5), $\alpha$ is unramified over
$K_P(\sqrt[p]{\pi}) \subset L_Q$ and hence unramified at $Q$.  Suppose
that char$(\kappa(P)) = p$.  Since $A$ is a complete regular local
ring of dimension 2 with maximal ideal $m = (\pi, \delta)$, we have
$A/P \simeq k[[\overline{\delta}]]$, where $\overline{\delta}$ is the
image of $\delta$ in $A/P$.  Since $k = k^p$, $[\kappa(P) :
\kappa(P)^p] = p$ and $\overline{\delta}$ is not a $p^{\rm th}$ power
in $\kappa(P)$.  Thus, by (1.2), $\alpha$ is split over $K_Q =
K_P(\sqrt[p]{\pi}, \sqrt[p]{\delta})$. In particular, $\alpha$ is
unramified at $Q$.  The case $P = (\delta)$ is similar.  \hfill $\Box$

\paragraph*{Proposition 2.7.} Let $A$ be a regular local ring of
dimension two with maximal ideal $m$ and $\pi \in m \setminus
m^2$. Let $K$ be the field of fractions of $A$ and $k$ the residue
field of $A$.  Suppose that char$(K) = 0$ and char$(k) = p>0 $.  Let
$\alpha \in \Brp(K)$. Suppose $\alpha$ is ramified on $A$ at most at
$(\pi)$.  Further assume that $\alpha$ is unramified over
$K_{\pi}(\sqrt[p]{\pi})$.  Then $\alpha$ is unramified over any
discrete valuation of $K(\sqrt[p]{\pi})$ dominating $A$.

\paragraph*{Proof.}  Let $B = A[\sqrt[p]{\pi}]$. Then $B$ is a regular
local ring with field of fractions $K(\sqrt[p]{\pi})$ and $B$ is
integral over $A$.  Every discrete valuation of $L = K(\sqrt[p]{\pi})$
dominating $A$ also dominates $B$. Thus, $B$ being regular, it is
enough to show that $\alpha$ is unramified on $B$.  Let $Q$ be a
height one prime ideal of $B$.  Since $B$ is integral over $A$, $Q
\cap A = P$ is a height one prime ideal. Suppose $P \neq (\pi)$. Since
$\alpha$ is unramified at $P$, $\alpha$ is unramified at $Q$.  Suppose
$P = (\pi)$.  Then by the assumption, $\alpha$ is unramified over $L_Q
= K_{\pi}(\sqrt{\pi})$ and hence unramified at $Q$ by (2.4).  \hfill
$\Box$

\vskip 5mm

\paragraph*{Lemma 2.8.}  Let $A$ be a one-dimensional Noetherian local
domain with field of fractions $K$ and residue field $k$.  Suppose
that char$(K) = p > 0$ and $[K : K^p] \leq p^{d+1}$ for some $d \geq
0$.  Then $[k : k^p] \leq p^d$.

\paragraph*{Proof.}  Let $B$ be the integral closure of $A$ in $K$ and
$m$ a maximal ideal of $B$. Then $B/m$ is a finite extension of $k$
and hence $[(B/m) : (B/m)^p] = [k : k^p]$. Thus, by replacing $A$ by
$B_m$, we assume that $A$ is a discrete valuation ring.  Let $\pi \in
A$ be a generator of the maximal ideal of $A$.  Let $u_1, \cdots , u_r
\in A^*$ be such that the images of $u_1, \cdots , u_r$ in $k$ are
$p$-independent.  Then, it is easy to see that $\pi, u_1, \cdots, u_r$
are $p$-independent in $K$.  \hfill $\Box$

\paragraph*{Theorem 2.9.} Let $X$ an excellent regular integral scheme
of dimension 2 and $K$ the function field of $X$.  Let $p$ be a prime
number. Suppose that char$(K) = 0$ and $K$ contains a primitive
$p^{\rm th}$ root of unity.  Suppose that for every codimension one
point $x$ of $X$ with char$(\kappa(x)) = p$, $[\kappa(x) :
\kappa(x)^p] = p$. Let $\alpha \in \Brp(K)$. Then there exist $f, g, h
\in K^*$ such that $\alpha \otimes K(\sqrt[p]{f}, \sqrt[p]{g},
\sqrt[p]{h})$ is unramified at every discrete valuation of
$K(\sqrt[p]{f}, \sqrt[p]{g}, \sqrt[p]{h})$ whose restriction to $K$ is
centered at a point of $X$.

\paragraph*{Proof.}  Let $P$ be a closed point of $X$. Suppose that
char$(\kappa(P)) = p$.  Let $\pi \in {\cal O}_{X, P}$ be a prime
dividing $p$.  Then $A = {\cal O}_{X, P}/(\pi)$ is a one-dimensional
Noetherian local domain. Let $\kappa(\pi)$ be the field of fractions
of $A$. Then $\kappa(\pi)$ is the residue field of a codimension one
point of $X$ and char$(\kappa(\pi)) = p$.  Thus, by the assumption,
$[\kappa(\pi) : \kappa(\pi)^p] = p$. Hence, by (2.8), $[\kappa(P) :
\kappa(P)^p] = 1$.  Let $X'$ be a blow-up of $X$ at finitely many
closed points.  Let $x' \in X'$ be a point of codimension 1.  If $x'$
is an exceptional curve, then $\kappa(x') = \kappa(P)(t)$ for some
closed point $P$ of $X$ and $t$ a variable.  In particular,
$[\kappa(x') : \kappa(x')^p] = p$.  Thus $X'$ also satisfies the
hypothesis of the theorem.

By blowing up $X$ at finitely many closed points ([Li]), we assume
that ram$(\alpha)$ and Supp$(p)$ are contained in $C + E$, where $C$
and $E$ are regular curves with $C$ intersecting $E$ transversally.
Let $D$ be an irreducible component in $C \cup E$. If char$(\kappa(D))
\neq p$, let $\pi_D$ be any parameter at $D$.  Then, by (2.5),
$\alpha$ is unramified over $K_D(\sqrt[p]{\pi_D})$. If
char$(\kappa(D)) = p$, then by the assumption, $[\kappa(D) :
\kappa(D)^p] = p$ and hence by (1.8) there exists a parameter $\pi_D
\in K_{D}^*$ at $D$ such that $\alpha$ is unramified over
$K_{D}(\sqrt[p]{\pi_D})$.  By the weak approximation, choose $f \in
K^*$ such that $f = \pi_D$ modulo $K_D^{*^p}$ for every irreducible
component $D$ in $C \cup E$. Then, we have $K_D(\sqrt[p]{f}) =
K_D(\sqrt[p]{\pi_D})$ and
$$div(f) = C + E + F$$
for some $F$ not containing any component of $C$ and $E$.  In
particular, $\alpha$ is unramified over $K_D(\sqrt[p]{f})$ for every
$D$ in $C \cup E$.

Let ${\cal P}$ be the finite set of closed points containing $C \cap
E$, $C \cap F$ and $E \cap F$ and at least one point from each
component of $C$, $E$, $F$. Let $A$ be the semi-local regular two
dimensional ring at the points in ${\cal P}$.  Let $\pi, \delta \in A$
be such that the divisor of $\pi$ on $A$ is $C$ and the divisor of
$\delta$ on $A$ is $E$. For $P \in {\cal P}$, let $m_P$ denote the
maximal ideal of $A$ at $P$.  If $P \in C \cap F$, $P \not\in E$, let
$\delta_P \in m_P$ be such that $\delta_P \not\in m_Q$ for all $Q \in
{\cal P}, Q \neq P$ and $\delta_P \not\in (\pi) + m_P^2$.  Similarly,
choose $\pi_Q$ for each $Q \in E \cap F$, $Q \not\in C$.

Let $$ g = \pi \prod_{Q \in (E \cap F) \setminus C} \pi_Q ~~{\rm
  and}~~ h = \delta \prod_{P \in (C \cap F) \setminus E} \delta_P.$$
Let $P \in {\cal P}$.  By the choice of $g$ and $h$, we have $m_P =
(g, h)$.  If $P \in C$, then $g$ defines $C$ at $P$ and if $P \in E$,
then $h$ defines $E$ at $P$.

We claim that $\alpha$ is unramified   over $L = K(\sqrt[p]{f},
\sqrt[p]{g}, \sqrt[p]{h})$.  Let $\nu$ be a
discrete valuation of $L$.  Let $R$ be the discrete valuation ring of
$L$ at $\nu$.   Then there is a point $x$
of $X$ such that $R$ dominates the local ring $A_x = {\cal O}_{X, x}$.

If $x$ is not on $C \cup E$, then $\alpha$ is unramified on $A_x$ andÄ
hence unramified at $\nu$.  Suppose  that $x \in C \cup E$.

Suppose that $x$ is a codimension one point.  Then $x$ corresponds to an
irreducible component of $C$ or $E$.  By the
choice of   $f$, $\alpha$ is unramified over $K_x(\sqrt[p]{f}) \subset
L_{\nu}$.
In particular $\alpha$ is unramified at $\nu$.

Suppose that $x$ is a closed point.  Suppose that $x \not\in  {\cal P}$.
 Since $x \in C \cup E$, $C$, $E$ are
regular curves on $X$,   $f \in
m_x\setminus m_x^2$ and   $\alpha$ is
ramified on ${\cal O}_{X, x}$ only along $(f)$. By the choice of
$f$ and by (2.7), $\alpha$ is unramified at $\nu$ restricted
to $K(\sqrt[p]{f})$. Since $K(\sqrt[p]{f}) \subset L$, $\alpha$ is
unramified at $\nu$.

Assume that $x \in {\cal P}$.  Then, by the choice of $g$ and $h$,   
the maximal ideal $m_x$ of ${\cal O}_{X, x}$ is generated by $
g$ and $h$. Further $\alpha$ is unramified on $A$ except at $(g)$ and $(h)$. 
Since $\nu$ is centered on ${\cal O}_{X, x}$, its
restriction to $K(\sqrt[p]{g}, \sqrt[p]{h})$ is also centered on
${\cal O}_{X, x}$. Hence, 
by (2.6), $\alpha$ is unramified at $\nu$ restricted to
$K(\sqrt[p]{g}, \sqrt[p]{h})$. Since  $K(\sqrt[p]{g}, \sqrt[p]{h})
\subset L$, $\alpha$ is unramified at $\nu$. 
\hfill $\Box$

\paragraph*{Corollary 2.10.}  Let $k$ be a $p$-adic field and $K$ a
function field of a curve over $k$. Then for every central simple
algebra $A$ over $K$,   index$(A)$  divides  (period$(A))^3$.

\paragraph*{Proof.}  Let $A$ be a central simple algebra over $K$ and
$\alpha$ its class in $Br(K)$.  Let $n$ be the period of $A$.  It is
enough to prove the result for a prime $n$ and assuming $K$ contains a
primitive $n^{\rm th}$ root of unity.  If $n$ is coprime to $p$, then
by ([S1]), period of $A$ divides the square of the index. Assume that
$n= p$.  Since $K$ is a field of fractions of a regular proper scheme
over the ring of integers in $k$ of dimension 2, by (2.9), there exist
$f, g, h \in K^*$ such that $\alpha $ is unramified over $K
(\sqrt[p]{f}, \sqrt[p]{g}, \sqrt[p]{h})$. Since the unramified Brauer
group of $K (\sqrt[p]{f}, \sqrt[p]{g}, \sqrt[p]{h})$ is zero ([G],
2.15 and 3.1), $\alpha$ is trivial over $K (\sqrt[p]{f}, \sqrt[p]{g},
\sqrt[p]{h})$. Hence the index of $\alpha$ divides $p^3$.  $\hfill
\Box$.

\paragraph*{Corollary 2.11.} Let $k$ be a number field and $K$ function
field of a curve over $K$. Then for every central simple algebra of
period $p$, there exist $f,g, h \in K^*$ such that $\alpha \otimes
K(\sqrt[p]{f}, \sqrt[p]{g}, \sqrt[p]{h})$ is unramified at every
discrete valuation of $K(\sqrt[p]{f}, \sqrt[p]{g}, \sqrt[p]{h})$.

\paragraph*{Proof.}  We assume that $K$ contains a primitive $p^{\rm
  th}$ root of unity.  Since $K$ is a field of fractions of a regular
proper scheme over the ring of integers in $k$ of dimension 2 and
every discrete valuation of $K$ is centered at a point of $X$, by
(2.9), there exist $f, g, h \in K^*$ such that $\alpha $ is unramified
over $K (\sqrt[p]{f}, \sqrt[p]{g}, \sqrt[p]{h})$ at every discrete
valuation of $K(\sqrt[p]{f}, \sqrt[p]{g}, \sqrt[p]{h})$ . \hfill
$\Box$

\section*{3. Period-index bounds}
In this section we ameliorate some of the results of [L2] on the
period-index problem for Brauer classes on arithmetic surfaces.  Let
$k$ be a number field, $\A$ the ring of ad\`eles of $k$ and $S$ the
scheme of integers of $k$.  We fix a regular proper surface $X$ with a
flat projective generically smooth morphism $X\to S$.  We will write
$C$ for the generic fiber of $X/S$, and we assume that $C$ is
geometrically connected.  Unlike in Section 7 of [L2], we do not
assume that $C$ admits a $k$-point; instead, we will use the following
Lemma to get $0$-cycles of degree $1$ on appropriate moduli spaces.

\paragraph*{Lemma 3.1}
 Let $\overline M$ be a proper smooth
connected $k$-scheme.  For every nonempty open subscheme
$M\subset\overline M$, we have that $M$  contains a $0$-cycle of
degree $1$ if and only if   $\overline M$ contains a $0$-cycle of degree $1$.

\paragraph*{Proof.}

It is enough to show that if $\overline M$ contains a point with
residue field a finite extension $L/k$ then $M$ contains a $0$-cycle
of degree equal to $[L:k]$.  Since $\overline M$ is smooth, it is easy
to see (by choosing generic parameters in the local ring at $Q$ and
taking a Zariski closure) that there is a smooth projective curve $D$
over $k$ with an $L$-rational point $R\in D$ and a finite morphism
$f:D\to\overline M$ such that $f(R)=Q$ and $U=f^{-1}(M)\neq\emptyset$.
Since any point of $U$ is ample, it is easy to see that there is a
divisor $E\subset U\subset D$ of degree $[L:k]$.  The proper
pushforward of $E$ to $\overline M$ gives a $0$-cycle of degree equal
to $[L:k]$ contained in $M$, as desired.  \hfill $\Box$
\vskip 5mm

Fix a prime number $\ell$ and a $\m_\ell$-gerbe $\ms X\to X$.  Write
$\ms C\to C$ for the restriction of $\ms X$ to $C$.  Let $\ms M:=\ms
M_{\ms C}(\ell)$ denote the stack of stable $\ms C$-twisted sheaves of
rank $\ell$ and trivial determinant (as described in [L1]).  We know
the following facts about this stack.

\begin{enumerate}
\item $\ms M$ is a $\m_\ell$-gerbe over a smooth quasi-projective
  variety $M$ admitting a locally factorial compactification $M\subset
  M^{ss}$ such that $$\operatorname{codim}(M^{ss}\setminus
  M,M^{ss})>2$$ and $$\Pic(M^{ss})=\Pic(M)=\Z.$$

\item the Brauer group $\Br(M\tensor\widebar k)$ is generated by the
  class associated to $\ms M\tensor \widebar k$, this class has period and
  index $\ell$, and the sequence
$$0\to\Br(k)\to\Br(M)\to\Br(M\tensor\widebar k)\to 0$$
is exact.  Note that since $\ms M$ is a $\m_\ell$-gerbe over $M$, its
Brauer class has period $\ell$ and thus the universal obstruction
provides a \emph{canonical\/} splitting of the above exact sequence.
\end{enumerate}

\paragraph*{Proposition 3.2}
If $k$ is totally imaginary then the Brauer-Manin set $M(\A)^{\Br}$ is non-empty.

\paragraph*{Proof.}  We first claim that for every place $v$ of $k$,
the category $\ms M(k_v)$ is non-empty.  First, suppose $v$ is
finite. Consider the base
change $\ms X\tensor_S\ms O_{k_v}\to X\tensor_S\ms O_{k_v}$.  This is
a $\m_\ell$-gerbe on a proper curve over a complete discrete valuation
ring with finite residue field.  Since the Brauer group of such a
scheme is trivial, we have an invertible $\ms C_{k_v}$-twisted sheaf
$\ms L$.  Writing $L:=(\ms L^{\tensor \ell})^\vee$, tensoring with
$\ms L$ defines an isomorphism
$$\ms M_{C_{k_v}}(\ell,L)\simto\ms M\tensor k_v.$$
Since stable vector bundles of rank $\ell$ and determinant $L$ exist
on any curve over an infinite field, the former has an object over
$k_v$, whence the latter does as well. When $v$ is infinite, the
completion is algebraically closed (as $k$ is assumed to be totally
imaginary), so we can use Tsen's theorem to trivialize the Brauer
class and similarly reduce to the existence of stable vector bundles
on curves over algebraically closed fields.

Write $(x_v)$ for the system of objects of $\ms M(k_v)$.  Projecting
to $M$ gives a point $(\widebar x_v)\in M(\A)$ with the property that
the pairing $(\widebar x_v)\cdot[\ms M\to M]$ is $0$ in $\Q/\Z$.  As
noted above, the Leray spectral sequence shows that $\Br(M)/\Br(k)$ is
generated by the class of $\ms M\to M$.  We conclude that the point
$(\widebar x_v)$ lies in $M(\A)^{\Br}$, as desired.  \hfill $\Box$

\paragraph*{Theorem 3.3}
Assuming the CT-conjecture (see page 2), if $k$ is totally imaginary 
any class $\alpha\in\Br(X)$ satisfies $\ind(\alpha) | \per(\alpha)^2$.

\paragraph*{Proof.}
By standard inductive arguments, we may assume that $\alpha$ has prime
index $\ell$.

We retain the notation from above.  Let $M\subset \overline M$ be any
smooth compactification of $M$; this is possible by Hironaka's
theorem, since $k$ has characteristic $0$.  By functoriality we have
that $M(\A)^{\Br}\subset\overline M(\A)^{\Br}$, hence by Propostion
3.2 $\overline M(\A)^{\Br}\neq\emptyset$.  The CT-conjecture applies
to show that $\overline M$ has a $0$-cycle of degree $1$.  Applying
Lemma 3.1, we have a $0$-cycle of degree $1$ in $M$.  Put another way,
there are two finite \'etale $k$-algebras $A_1$ and $A_2$ such that
$[A_1:k]$ is relatively prime to $[A_2:k]$ and such that
$M(A_1)\neq\emptyset$ and $M(A_2)\neq\emptyset$.

Given a finite \'etale closed subscheme $Z\in M$, the Brauer class of
$\ms M\times_M Z$ lies in $\Br(Z)[\ell]$ and thus has index dividing
$\ell$ (as the residue fields of the Artinian scheme $Z$ are global
fields).  It follows that there is a finite flat $Z$-scheme $Y\to Z$
of degree $\ell$ such that $\ms M(Y)\neq\emptyset$.  Applying this to
$A_1$ and $A_2$ above, we find that for $i=1,2$ there is a finite flat
$A_i$-algebra $B_i$ of degree $\ell$ and a locally free $\ms
C_{B_i}$-twisted sheaf of rank $\ell$.  Putting these together and
arguing as in Section 4.1.1 of [L3], we see that the index of $\alpha$
divides $\ell^2[A_i:k]$ for $i=1,2$.  Since $[A_1:k]$ and $[A_2:k]$
are relatively prime, we see that $\ind(\alpha) |\ell^2$, as desired.
\hfill $\Box$

\section*{4. $u$-invariant of the function field a curve over a number
field}

\paragraph*{Theorem 4.1.} Assume that the CT-conjecture holds.  Let
$k$ be a totally imaginary number field and $K$ the function field of a curve over $k$.
Then for any central simple algebra $A$ over $K$, ind$(A)$ divides
period$(A)^5$.

\paragraph*{Proof.}  Let $A$ be a central simple algebra over $K$.  By
a standard inductive argument, it is enough to consider the case
period$(A) = p$ a prime and $K$ contains a primitive $p^{\rm th}$ root
of unity.  By (2.11), there exist $f, g, h \in K^*$ such that $A
\otimes_K K(\sqrt[p]{f}, \sqrt[p]{g}, \sqrt[p]{h})$ is unramified at
every discrete valuation of $K(\sqrt[p]{f}, \sqrt[p]{g},
\sqrt[p]{h})$.  By (3.3), index$(A\otimes K(\sqrt[p]{f}, \sqrt[p]{g},
\sqrt[p]{h}))$ divides $p^2$. Hence index$(A)$ divides $p^5$. \hfill
$\Box$

\paragraph*{Corollary 4.2.} Assume that the CT-conjecture holds.  Let
$k$ be a number field and $K$ the function field of a curve over $k$.
Then for any central simple algebra $A$ over $K$, ind$(A)$ divides
period$(A)^6$.

\paragraph*{Proof.}  By (4.1),  ind$(A \otimes_kk(\sqrt{-1}))$ divides $2^5$. 
Hence  ind$(A)$ divides $2^6$. \hfill $\Box$

\paragraph*{Theorem 4.3.} Assume that the CT-conjecture holds.  Let
$k$ be a number field and $K$ the function field of a curve over
$k$. Then there exists an integer $N_2$ (which does not depend on $K$
or $k$) such that every element in $H^2(K, \mu_2)$ is a sum of at most
$N_2$ symbols.

\paragraph*{Proof.}  Let $\alpha \in H^2(K, \mu_2)$. By (4.2),
index$(\alpha)$ divides $2^6$.  Let ${\cal A}$ be the generic division
algebra of degree $2^6$ with center $Z$ and $X_{{\cal A}^2}$ be the
Severi-Brauer variety of ${\cal A} \otimes {\cal A}$ over $Z$. Let
$\alpha$ be represented by a central simple algebra $A$ over $K$ of
degree $2^6$.  Then there is a specialisation from the function field
$Z(X_{{\cal A}^2})$ to $K$ which specialises ${\cal A}$ to $A$. Since
${\cal A}$ is 2-torsion element in the Brauer group of $Z(X_{{\cal
    A}^2})$, ${\cal A}$ is a product of $N_2$ quaternion algebras over
$Z(X_{{\cal A}^2})$ ([M]). In particular $A$ is a product of $N_2$
quaternion algebras over $K$.  Hence $\alpha$ is a sum of at most
$N_2$ symbols. \hfill $\Box$.

\paragraph*{Theorem 4.4.} Assume that the  CT-conjecture
 holds.  Let $k$ be a  number field and $K$ the function field
of a curve over $k$.  Then every element in $H^3(K, \mu_2)$ is a  sum
of $N_2$ symbols.

\paragraph*{Proof.} Let $\beta \in H^3(K, \mu_2)$. Then by ([Su]),
there exists $f \in K^*$ such that $\alpha$ is zero over
$K(\sqrt{f})$. Hence $\beta = (f) \cdot \alpha$ for some $\alpha \in
H^2(K, \mu_2)$ ([A], 4.6). By (4.3), $\alpha = (a_1) \cdot (b_1) +
\cdots + (a_{N_2}) \cdot (b_{N_2})$. Thus $\beta = (f) \cdot \alpha =
(f)\cdot (a_1) \cdot (b_1) + \cdots + (f) \cdot (a_{N_2}) \cdot
(b_{N_2})$. \hfill $\Box$

\paragraph*{Theorem 4.5} Assume that the  CT-conjecture  holds. 
Let $k$ be a totally imaginary  number field and $K$ the function field
of a curve over $k$. Then $u(K)$ is finite.

\paragraph*{Proof.}  By (4.3), there exists $N_2$ such that every
element in $H^2(K, \mu_2)$ is a sum of at most $N_2$ symbols.  By
(4.4), every element in $H^3(K, \mu_2)$ is a sum of at most $N_2$
symbols.  Since $k$ is a totally imaginary number field, $H^4(K,
\mu_2) = 0$. Hence $u(K)$ is finite (cf. [PS]). \hfill $\Box$

\section*{References}

\begin{enumerate}

\item[{[A]}] J.\ K.\ Arason, \emph{Cohomologische Invarianten
    quadratischer Formen}, J.Algebra {\bf 36} (1975), 448--491.

\item[{[AG]}] M.\ Auslander, O.\ Goldman, \emph{Maximal orders}, 
Trans.\ Amer.\ Math.\ Soc.\ {\bf 97} (1960), 367--409.

\item[{[APS]}] Asher Auel, R. Parimala and V. Suresh, \emph{Quadric
    surface bundles over surfaces}, arXiv:1207.4105v1

\item[{[B]}] E. Brussel, \emph{On Saltman's p-adic curves papers},
  Quadratic forms, linear algebraic groups, and cohomology, 13-39,
  Dev. Math., {\bf 18}, Springer, New York, 2010.

\item[{[CT1]}] J.-L. Colliot-Th\'el\`ene \emph{Conjectures de type
    local-global sur image des groupes de Chow dans la cohomologie
    \'etale}, Algebraic K -theory (Seattle, WA, 1997), 1-12,
  Proc. Sympos. Pure Math., {\bf 67}, Amer. Math. Soc., Providence,
  RI, 1999.

\item[{[CT2]}] J.-L. Colliot-Th\'el\`ene \emph{Cohomologie des corps
    valu\'es henséliens, d'apr\`es K. Kato et S. Bloch}, Algebraic
  $K$-theory and its applications (Trieste, 1997), ed. H. Bass,
  A. O. Kuku, C. Pedrini, World Scientific Publishing, 1999,
  p. 120-163.

\item[{[CTS]}] J.-L. Colliot-Th\'el\`ene et J.-J.  Sansuc,
  \emph{Fibr\'es quadratiques et composantes connexes r\'eelles},
  Math. Ann.  {\bf 244} (1979), 105-134.

\item[{[G]}] A. Grothendieck, \emph{Le groupe de Brauer III}, Dix
  exposés sur la cohomologie des schémas, Amsterdam, North-Holland,
  Amsterdam (1968), 88-188.

\item[{[K]}] K. Kato, \emph{Galois cohomology of complete discrete
    valuation fields}, in Algebraic K-Theory, Springer L. N. M. {\bf
    967} (1982), 215-238.

\item[{[L1]}] M. Lieblich, \emph{Moduli of twisted sheaves}, Duke
  Math. J. {\bf 138} (2007), 23-118.
    
\item[{[L2]}] M. Lieblich, \emph{Period and index in the Brauer group
    of an arithmetic surface}, with an appendix by Daniel Krashen,
  J. Reine Angew. Math {\bf 659} (2011), 1-41.

\item[{[L3]}] M. Lieblich, \emph{Twisted sheaves and the period-index
    problem}, Compos. Math. {\bf 144} (2008), no.\ 1, 1-31.
    
\item[{[Li]}] J. Lipman, \emph{Desingularization of two-dimensional
    schemes}, Ann. Math. {\bf 107} (1978), 151Ð207

\item[{[M]}] A.S. Merkurjev, \emph{On the norm residue symbol of
    degree 2}, Dokl. Akad. Nauk SSSR {\bf 261} (1981), 542-547.

\item[{[MS]}] A.S. Merkurjev and A. A. Suslin,  \emph{K-cohomology  of 
Severi--Brauer varieties and the norm- residue homomorphism}, Izv. Akad. 
Nauk SSSR, Ser. Mat., {\bf 46}, No. 5,  (1982) 1011-1061.

\item[{[S1]}] D.J. Saltman, \emph{Division Algebras over $p$-adic
    curves}, J. Ramanujan Math. Soc. {\bf 12} (1997), 25-47.

\item[{[S2]}] D.J. Saltman, \emph{Cyclic Algebras over $p$-adic
    curves}, J. Algebra {\bf 314} (2007) 817-843.

\item[{[S3]}] D.J. Saltman, \emph{ Division algebras over surfaces},
  J. Algebra {\bf 320} (2008), 1543-1585.

\item[{[S4]}] D.J. Saltman, \emph{ Bad characteristic}, private notes.

\item[{[PS]}] R. Parimala and V. Suresh, \emph{On the length of a
    quadratic form}, Algebra and number theory, 147-157, Hindustan
  Book Agency, Delhi, 2005.

\item[{[Su]}] V. Suresh, \emph{Galois cohomology in degree 3 of
    function fields of curves over number fields}, J. Number Theory
  {\bf 107} (2004), no. 1, 80-94.

\end{enumerate} 

\noindent
Department of Mathematics,  University of Washington, Box 354350 Seattle, WA 98195 

\vskip 3mm
\noindent
Department of Mathematics \& Computer Science, Emory University, 400 Dowman Drive NE W401, Atlanta, GA 30322

\vskip 3mm

\noindent
E-mail addresses:  lieblich@math.washington.edu,  parimala@mathcs.emory.edu,  suresh@mathcs.emory.edu

\end{document}